\documentclass[12pt,reqno]{amsart}

\usepackage{amsthm,amsmath,amssymb,amsfonts,amscd,verbatim,latexsym,xypic}
\newtheorem{Theorem}{Theorem}[section]
\newtheorem{Lemma}[Theorem]{Lemma}
\newtheorem{Proposition}[Theorem]{Proposition}
\newtheorem{Corollary}[Theorem]{Corollary}
\newtheorem{Remark}[Theorem]{Remark}

\newcommand{\be} {\begin{equation}}
\newcommand{\ee} {\end{equation}}
\newcommand{\bea} {\begin{eqnarray}}
\newcommand{\eea} {\end{eqnarray}}
\newcommand{\lp}  {\left(}
\newcommand{\rp}  {\right)}

\newcommand{\Complex}{\mathbf C}
\newcommand{\Proj}{\text{Proj} \, }
\renewcommand{\P}{\mathbb P}
\renewcommand{\O}{\mathcal O}
\newcommand{\U}{\mathcal U}
\newcommand{\E}{\mathcal E}
\newcommand{\F}{\mathbb F}
\newcommand{\G}{\mathbb G}
\newcommand{\Sz}{\mathcal S}
\newcommand{\Be}{\mathbb B}
\newcommand{\ux}{\mathbf x}
\newcommand{\uy}{\mathbf y}
\newcommand{\uz}{\mathbf z}

\newcommand{\im}{\text{im}}
\newcommand{\He}{\mathbb H}
\newcommand{\T}{\mathbb T}

\newcommand{\rank}{\text{rank}\,}
\newcommand{\Hom}{\text{Hom}}
\newcommand{\ord}{\text{ord}\,}
\newcommand{\codim}{\text{codim}\;}

\renewcommand{\theta}{\vartheta}
\newcommand{\R}{\mathfrak R}
\newcommand{\ga}{\mathfrak I}
\newcommand{\gb}{\mathfrak b}
\newcommand{\gP}{\mathfrak P}
\newcommand{\gq}{\mathfrak q}
\newcommand{\constant}{\Bbbk \;}
\newcommand{\A}{\mathcal A}
\newcommand{\zA}{M}
\newcommand{\zB}{N}
\newcommand{\bA}{\mathbb A}
\newcommand{\sC}{\mathcal C}
\newcommand{\tGamma}{{\widetilde \Gamma}}
\newcommand{\tPsi}{{\widetilde \Psi}}
\newcommand{\tF}{\widetilde F}
\newcommand{\tG}{\widetilde G}
\newcommand{\Mor}{\mathcal M}

\newcommand{\DD}{\mathcal  D}
\newcommand{\bone}{\mathbf 1}
\newcommand{\ra}{\rightarrow}

\newcommand{\lra}{\longrightarrow}
\newcommand{\la}{\leftarrow}

\newcommand{\demo}{\noindent {\sc Proof.}\;}

\begin{document}
\title{On the Jacobian ideal of the binary discriminant}
\author[D'Andrea and Chipalkatti]{Carlos D'Andrea and Jaydeep Chipalkatti} 
\maketitle 
\vspace{-0.6cm} 
\centerline{(with an appendix by {\sc Abdelmalek Abdesselam})}

\vspace{1.5cm} 

\parbox{12.5cm}{\small 
{\sc Abstract.} 
Let $\Delta$ denote the discriminant of the generic binary $d$-ic. 
We show that for $d \ge 3$, the Jacobian ideal of $\Delta$ is perfect of height 
$2$. Moreover we describe its $SL_2$-equivariant minimal resolution and 
the associated differential equations satisfied by $\Delta$. A similar 
result is proved for the resultant of two forms of orders $d,e$ whenever $d \ge e-1$. 
If $\Phi_n$ denotes the locus of binary forms with total root multiplicity $\ge d-n$, 
then we show that the ideal of $\Phi_n$ is also perfect, and we construct a covariant 
which characterizes this locus. We also explain the role of the Morley 
form in the determinantal formula for the resultant. This relies upon 
a calculation which is done in the appendix by A.~Abdesselam.}

\bigskip 

\parbox{12cm}{\small 
} 

\medskip 

\section{Introduction} 
\subsection{} Let 
\[ \F = a_0 \, x_1^d + \dots + 
\binom{d}{i} \, a_i \, x_1^{d-i} \, x_2^i + \dots + a_d \, x_2^d, \] 
denote the generic binary form of order $d$ in the variables $x_1,x_2$. 
Its 
discriminant $\Delta = \Delta(a_0,\dots,a_d)$ is a homogeneous polynomial
 with the following property: given 
$\alpha_0,\dots,\alpha_d \in \Complex$, the form 
$F_\alpha = \sum\limits_{i=0}^d \, \binom{d}{i} \, \alpha_i \, x_1^{d-i} x_2^i$ 
is divisible by the square of a linear form iff 
$\Delta(\alpha_0,\dots,\alpha_d)=0$. 
Let $R$ denote the polynomial ring $\Complex[a_0,\dots,a_d]$, and 
let 
\[ J = (\frac{\partial \Delta}{\partial a_0}, \dots, 
\frac{\partial \Delta}{\partial a_d}) \subseteq R, \] 
denote the Jacobian ideal of $\Delta$. 

Our main result (in \S\ref{section.J_Delta}) is that 
$J$ is a {\sl perfect} ideal of height $2$ 
for $d \ge 3$, with graded minimal resolution 
\begin{equation} 
 0 \la R/J \la R \la R(3 - 2d)^{d+1} \la R(2-2d)^{3} \oplus R(1-2d)^{d-3} \la 0. 
\label{res1} \end{equation}

\subsection{} \label{section.defnXlambda}
To put this statement into a geometric context, identify the form 
$F_\alpha$ (distinguished up to a scalar) with the point 
$[\alpha_0, \dots, \alpha_d]$ in the projective space $\P^d$. 
We recall the notion of a Coincident Root locus introduced in~\cite{ego3}. Let 
\[ \lambda = (\lambda_1, \lambda_2,\dots, \lambda_n) \] 
be a partition of $d$ into $n$ parts. Now the CR locus associated to $\lambda$ is 
defined to be 
\[ X_\lambda = \{ F \in \P^d: F = \prod\limits_{i=1}^n \, l_i^{\lambda_i} \; \; 
\text{for some linear forms $\l_i$} \},   \] 
which is an irreducible projective subvariety of dimension $n$. Given two 
partitions $\lambda$ and $\mu$, we have $X_\lambda \subseteq X_\mu$ iff $\mu$ is a 
refinement of $\lambda$. Now $X_{(2,1^{d-2})}$ is the hypersurface $\{\Delta =0\}$, 
and the closed subscheme $Z = \Proj (R/J)$ is supported on its singular locus. 
By~\cite[Theorem 5.4]{ego3}, the latter is equal to the union 
$X_\tau \cup X_\delta$, 
where 
\[ \tau = (3,\underbrace{1,\dots,1}_{d-3}) \quad \text{and} \quad 
\delta = (2,2,\underbrace{1,\dots,1}_{d-4}). \] 
The result above implies that $Z$ is an arithmetically Cohen-Macaulay scheme. 
In Proposition~\ref{prop.multiplicity} we show that $Z$ has multiplicities
$2$ and $1$ along $X_\tau$ and $X_{\delta}$ respectively. 

\subsection{} 
The ideas in \S\ref{section.J_Delta} are based on the `Cayley method'  
as explained in~\cite[Ch.~2]{GKZ}. In \S\ref{section.J_Res} we give a 
pr{\'e}cis of this method in the context of binary resultants, and then deduce the 
following theorem:  let $\R$ denote the resultant of generic binary 
forms $\F,\G$ of orders $d,e$. If $d \ge e-1$, then the $\F$-Jacobian ideal of 
$\R$ (i.e., the ideal of partial derivatives of $\R$ with respect to the coefficients of 
$\F$) is perfect. The Cayley method involves constructing a morphism of vector 
bundles whose determinant is the resultant. The most interesting ingredient in this morphism 
is the so-called Morley form $\Mor$, which encodes the $d_2$-differential 
of a spectral sequence. Although \emph{a priori} the differential is only well-defined 
modulo coboundaries, it admits a unique equivariant lifting to 
a morphism from binary forms of order $e-2$ to those of order $d$. 
This is explained in~\S\ref{proof.prop.morley} -- \ref{section.FGrGFr}, 
modulo a calculation which is provided in the appendix by A.~Abdesselam. 
The reader may also consult~\cite[\S3.11]{Jouanolou} for a very general treatment of 
multivariate Morley forms. 

\subsection{} 
In a slightly different direction, define 
$\Phi_n = \bigcup\limits_\lambda \,  X_\lambda$, 
where the union is quantified over all partitions $\lambda$ having $n$ parts. 
E.g., for $d=6$ and $n=3$, 
\[ \Phi_3 = X_{(4,1,1)} \cup X_{(3,2,1)} \cup X_{(2,2,2)}. \] 
Let $I_n \subseteq R$ denote the ideal of $\Phi_n$. In~\S\ref{section.acm} 
we show that $I_n$ is a determinantal ideal which admits an Eagon-Northcott 
resolution, in particular it is perfect.

\subsection{} 
Note that the group $SL_2 \, \Complex$ acts on $\P^d$, namely 
the element 
\[ g = \left(\begin{array}{rr} p & r \\ q & s \end{array} \right) 
\in SL_2 \, \Complex, \] 
sends $\sum\limits_i \, \binom{d}{i} \, \alpha_i \, x_1^{d-i} \, x_2^i$ to 
$\sum\limits_i \, \binom{d}{i} \, \alpha_i \, (p \, x_1+q \, x_2)^{d-i} \, 
(r \, x_1+s \, x_2)^i$. 
All the varieties defined above inherit this action, in particular 
the ideals $I_n,J$ and the Betti modules in their free resolutions are 
$SL_2$-representations. This equivariance is respected in all of our subsequent 
constructions. The first syzygy modules occuring in the resolution of $J$ encode the 
invariant differential equations satisfied by $\Delta$ (and similarly for $\R$). 
We write down these equations explicitly using transvectants. 
The reader is referred to~\cite[Lecture 11]{FH} 
and~\cite[\S 4.2]{Sturmfels} for basic representation theory of 
$SL_2$. We will use~\cite{Glenn} and~\cite{GrYo} as standard references for 
classical invariant theory and symbolic calculus; 
more recent accounts of this subject may be found in~\cite{CD_inv, Dolgachev1,KR, Olver}.

\medskip 

\noindent {\small 
{\sc Acknowledgements.} We thank Bernd Sturmfels for initiating the 
collaboration which led to this paper. We arrived at many of the 
results in this paper by extensive calculations in Macaulay-2, and it is 
a pleasure to thank its authors Dan Grayson and Mike Stillman. The second 
author was supported by NSERC while this work was in progress.} 

\section{Preliminaries} 
Let $V$ be a two-dimensional vector space over $\Complex$ 
with basis $\ux = \{x_1,x_2\}$. Then $\text{Sym}^m \, V = S_m \, V$ is the 
$(m+1)$-dimensional space of binary forms of order $m$ in $\ux$. 
The $\{S_m \, V : m \ge 0\}$ are a complete set of irreducible $SL(V)$-representations. 
We will omit the $V$ if no confusion is likely, thus $S_m(S_n)$ stands 
for the plethysm representation $\text{Sym}^m \, (\text{Sym}^n \, V)$ etc. 

\subsection{Transvectants} \label{section.trans} 
Given integers $m,n \ge 0$, we have a decomposition of 
$SL_2$-representations 
\begin{equation} 
S_m \otimes S_n \simeq \bigoplus\limits_{r=0}^{\min\{m,n\}} \, 
S_{m+n-2r}.   \label{Clebsch-Gordan} \end{equation}
Let $A,B$ denote binary forms of respective orders $m,n$. The $r$-th transvectant 
of $A$ with $B$, written $(A,B)_r$, is defined to be the image of 
$A \otimes B$ via the projection map 
\[ S_m \otimes S_n \lra S_{m+n-2r} \, .  \] 
It is given by the formula 
\begin{equation} (A,B)_r = \frac{(m-r)! \, (n-r)!}{m! \, n!} \, 
\sum\limits_{i=0}^r \, (-1)^i \binom{r}{i} \, 
\frac{\partial^r A}{\partial x_1^{r-i} \, \partial x_2^i} \, 
\frac{\partial^r B}{\partial x_1^i \, \partial x_2^{r-i}} 
\label{trans.formula} \end{equation} 
By convention $(A,B)_r = 0$ if $r > \min \, \{m,n\}$. (Some authors choose the 
scaling factor differently, cf.~\cite[Ch.~5]{Olver}.) 
Each $S_m$ is isomorphic to its dual representation 
$S_m^* = \text{Hom}(S_m,S_0)$ by the map 
which sends $A \in S_m$ to the functional $B \lra (A,B)_m$.  
Two forms $A,B \in S_m$ are said to be {\sl apolar} to each other if 
$(A,B)_m=0$. In some of the examples below quite a few complicated transvectants had to 
be calculated; to this end we programmed formula~(\ref{trans.formula}) in {\sc Maple}. 
If two forms are symbolically expressed, a useful general procedure for calculating 
their transvectants is given in~\cite[\S 3.2.5]{Glenn} (also see~\cite[\S 49]{GrYo}). 

\subsection{} We identify the generic binary $d$-ic 
$\F = \sum\limits_{i=0}^d \, \binom{d}{i} \, a_i \, x_1^{d-i} \, x_2^i$ 
with the natural trace form in 
$S_d \, \otimes \, S_d^*$. Using the self-duality above, this amounts to the identification 
of $a_i \in S_d^*$ with $\frac{1}{d!} \, x_2^{d-i} \, (-x_1)^i$. 
Let $R$ be the symmetric algebra 
\[ \bigoplus\limits_{m \ge 0} \, S_m(S_d^*) = \bigoplus\limits_{m \ge 0} \, R_m = 
\Complex \, [a_0,\dots,a_d], \] 
and $\P^d = \P \, S_d = \Proj \, R$. Generally $F,G,\dots$  will denote specific 
binary forms, as opposed to generic forms $\F,\G, \dots$. 

\subsection{} \label{section.cov}
A {\sl covariant} of degree-order $(m,q)$ of binary $d$-ics is 
by definition a trivial summand in the representation 
$S_q \otimes R_m$ (cf.~\cite[\S 11 et seq.]{GrYo}). An invariant is a covariant of order $0$. 
The most frequently appearing covariants are the Hessian $\He = (\F,\F)_2$, and 
the cubicovariant $\T = (\F,\He)_1$, of degree-orders $(2,2d-4)$ and 
$(3,3d-6)$ respectively. The discriminant $\Delta$ is 
an invariant of order $2(d-1)$. If $I(a_0,\dots,a_d)$ is an invariant of degree $m$, then its 
{\sl evectant} is defined to be 
\[ 
\E_I = \frac{(-1)^d}{m} \sum\limits_{i=0}^d \, \frac{\partial I}{\partial a_i} \, 
x_2^{d-i} \, (-x_1)^i. \] 
It is a covariant of degree-order $(m-1,d)$. The scaling factor is so 
chosen that we have an identity $(\E_I,\F)_d= I$. 

\subsection{} \label{deg.Xlambda}
The degree of the CR locus $X_\lambda$ is given by a formula due to Hilbert~\cite{Hilbert2}. 
Let $e_r$ denote the number of parts in $\lambda$ equal to $r$, thus 
$\sum\limits_{r \ge 1} e_r = n$ and $\sum r \, e_r =d$. Then 
$\deg X_\lambda = \frac{n!}{\prod\limits_r \, (e_r!)} \, \prod\limits_{i=1}^n \lambda_i$.  
For instance, 
$\deg X_{(3^2,2,1^3)} = \frac{6!}{2! \, 1! \, 3!} \; {3^2 \times 2 \times 1^3} = 
1080$. 

\section{The binary discriminant} \label{section.J_Delta} 
Throughout this paper, we will regard $\Delta$ and $\R$ as well-defined only up to a 
multiplicative constant. 

For a binary $d$-ic $F$, we define its Bezoutiant $\Be_F$ as follows: introduce new variables 
$\uy = (y_1, y_2)$, and write $G$ for the form obtained by 
substituting $y_1,y_2$ for $x_1,x_2$ in $F$. Then 
\[ \Be_F = (\frac{\partial F}{\partial x_1} \frac{\partial G}{\partial y_2} - 
\frac{\partial G}{\partial y_1} \frac{\partial F}{\partial x_2})/(x_1 \, y_2 - x_2 \, y_1),  
\] 
which is a form of order $(d-2,d-2)$ in $\ux,\uy$. 
Henceforth we will assume $d \ge 4$ (but see \S\ref{d_2_or_3}). In the sequel, 
$\Bbbk$ will stand for a nonzero rational constant which need not be precisely specified. 
Define a map 
\[ \beta_F: S_{d-4} \lra S_d, \] 
by sending $A \in S_{d-4}$ to $[(A,\Be_F)_{d-4}]_{\uy=\ux}$. This is 
interpreted as follows: 
take the $(d-4)$-th transvectant of $A$ with $\Be_F$ with 
respect to the $\ux$ variables, which gives an $\ux \, \uy$-form of order 
$(2,d-2)$. By substituting $\ux$ for $\uy$ we get an $\ux$-form of order $d$. 
Define another morphism \[ 
\gamma_F: S_2 \lra S_d, \quad A \lra (A,F)_1,  \] 
and finally let 
\[ \bone_F: S_0 \lra S_d, \quad 1 \lra F. \] 
Note that $\beta_F$ is quadratic in the coefficients of $F$, whereas 
$\gamma_F,\bone_F$ are linear. Now consider the morphism 
\[ \underbrace{\, \beta_F \oplus \gamma_F \oplus \bone_F}_{h_F}: 
S_{d-4} \oplus S_2 \oplus S_0 \lra S_d. \] 

\begin{Proposition} \sl 
We have an equality 
\[ \det \, h_\F  = \Delta_\F \] 
up to a nonzero scalar. 
\end{Proposition} 

\demo Let $D_\F = \det h_\F$. It is an invariant of 
degree $2(d-3) + 3 + 1 = 2(d-1)$, which is the same as $\deg \Delta_\F$. 
We will show that 
(1) $D_\F$ vanishes whenever $F$ has a repeated linear factor, and 
(2) $D_\F$ is not identically zero. This will imply that $D_\F = \Delta_\F$ 
(up to a scalar). 

As to (1), after a change of variables we may assume that $x_1^2$ divides $F$. 
Then $x_1 \, y_1$ divides $\Be_F$, and hence $x_1$ divides each form 
in $\im(\beta_F)$. 
Similarly, $x_1$ divides each form in $\im(\gamma_F)$ and $\im(\bone_F)$, hence 
$h_F$ is not surjective and $D_F = 0$. 
Now assume $F = x_1^d + x_2^d$, then 
\[ \Be_F = d^2 \, \sum\limits_{i=0}^{d-1} \, 
(x_1 \, y_2)^{d-2-i} \, (x_2 \, y_1)^i.  \] 
By a direct calculation, $\beta_F(x_1^{d-k-4}x_2^k) = 
\constant x_1^{d-k-2} \, x_2^{k+2}$, 
hence $\im(\beta_F) = \text{Span} \, \{ x_1^{d-i} x_2^i: 2 \le i \le d-2\}$. Since 
\[ \gamma_F(x_1^2) = \constant x_1 \, x_2^{d-1}, \quad 
  \gamma_F(x_1 \, x_2) = \constant (x_1^d - x_2^d), \quad 
  \gamma_F(x_2^2) = \constant x_1^{d-1} \, x_2,\] 
we deduce that $h_F$ is surjective. This shows (2) and completes the proof. 
\qed 

\smallskip 

A similar calculation shows that if $F = x_1^2 \, (x_1^{d-2} + x_2^{d-2})$, then 
$\im(h_\F) = \text{Span} \, \{x_1^{d-i} \, x_2^i: 0 \le i \le d-1\}$. Hence 
$h_F$ has rank $d$ for a general $F \in X_{(2,1^{d-2})}$. 
Let $\E_\Delta$ be the evectant of $\Delta$ (see \S \ref{section.cov}), 
and define the map 
\[ e_F: S_d \lra S_0, \quad A \lra (A,\E_\Delta)_d. \] 

\begin{Lemma} \sl 
The composites 
\[ e_F \circ \beta_F: S_{d-4} \lra S_0, \quad 
e_F \circ \gamma_F: S_2 \lra S_0 \] 
are zero. 
\end{Lemma} 
\demo Since $e_F \circ \beta_F$ is of degree $(2d-1)$ in the coefficients of $F$, 
it corresponds to an $SL_2$-equivariant map $S_{d-4} \lra R_{2d-1}$. 
Said differently, there exists 
a covariant $C$ of $d$-ics of degree-order $(2d-1,d-4)$ such that 
$e_F \circ \beta_F(A) = (A,C)_{d-4}$. Similarly, there is a $C'$ of degree-order 
$(2d-2,2)$ such that $e_F \circ \gamma_F(A) = (A,C')_2$. 

We will show that 
if $F \in X_{(2,1^{d-2})}$, then $e_F \circ \beta_F = e_F \circ \gamma_F = 0$. 
This will imply that each coefficient of $C$ or $C'$ vanishes on $X_{(2,1^{d-2})}$, 
and hence must be divisible by $\Delta_\F$. The quotients 
$C/\Delta_\F, C'/\Delta_\F$ are of degree-orders $(1,d-4)$ and $(0,2)$ respectively. 
Since there are no such nonzero covariants, $C$ and $C'$ must be zero. 

Let $x_1^2$ be a factor of $F$. By~\cite[Ch.~12, formula (1.28)]{GKZ} 
(also see~\cite[Art.~96]{Salmon1}), 
we have $\E_\Delta = \constant x_1^d$. Any form $B$ in the image of 
$\beta_F$ or $\gamma_F$ is divisible by $x_1$, 
hence $(B,\E_\Delta)_d = (B, x_1^d)_d = 0$. This completes the proof. \qed 

\subsection{} 
Now consider the map 
\[ \beta_\F \oplus \gamma_\F: S_{d-4} \oplus S_2 \lra S_d, \] or 
what is the same, the corresponding map of graded $R$-modules 
\begin{equation} R(-2) \otimes S_{d-4} \oplus R(-1) \otimes S_2 \lra R \otimes S_d.
\label{map.Rmodules} \end{equation}
Let $M$ denote its $d \times (d+1)$ matrix with respect to the natural monomial bases. 

\begin{Lemma} \sl 
The ideal of maximal minors of $M$ equals $J$ (the Jacobian ideal of $\Delta$). 
\end{Lemma} 
\demo 
Let $W$ denote the image of $1$ via the map 
\[ \wedge^d \, (\beta_\F \oplus \gamma_\F): \Complex \lra \wedge^d S_d \simeq S_d. \] 
By construction $W$ is a covariant of degree-order $(2d-3,d)$ whose coefficients are 
exactly the maximal minors. Let $\{ A_1,\dots,A_d\}$ span 
$\im(\beta_\F \oplus \gamma_\F)$. 
On the one hand, $W$ is the Wronskian of the $A_i$, hence it is (up to scalar) the 
unique $d$-ic which is apolar to all the $A_i$ (see \cite[Appendix II]{GrYo}). 
On the other hand, $(A_i,\E_\Delta)_d =0$ by the lemma above. 
Hence $W = \constant \E_\Delta$. \qed 

\medskip 

The subvariety of $\P^d$ defined by $J$ is codimension $2$, 
hence the Eagon-Northcott complex (or what is the same in this case, the Hilbert-Burch 
complex) of the map~(\ref{map.Rmodules}) resolves $J$ (see~\cite[Ch.~16 F]{BrunsVetter}). 
We have proved the following: 
\begin{Theorem} \label{theorem.res.J} 
\sl The ideal $J$ is perfect of height $2$ with 
$SL_2$-equivariant minimal resolution 
\[ \begin{aligned} 
0 \la R/J \la R & \la R(3-2d) \otimes S_d \\ 
& \la R(2-2d) \otimes S_2 \oplus R(1-2d) \otimes S_{d-4} \la 0. \qquad \qed 
\end{aligned}  \]
\end{Theorem}

\subsection{} The first syzygy modules $S_2, S_{d-4}$ correspond to systems of 
$SL_2$-equivariant differential equations for $\Delta$, we proceed to make these equations explicit. 
For all $A \in S_2$, we have $((A,\F)_1,\E_\Delta)_d =0$. 
Using classical symbolic calculus (see~\cite[Ch.~I]{GrYo}), let 
\[ A = \alpha_\ux^2, \quad \F = f_\ux^d, \quad \E_\Delta = e_\ux^d. \] 
Then $(A,\F)_1 = (\alpha \, f) \, \alpha_\ux \, f_\ux^{d-1}$, and 
\[ \begin{aligned} 
{} & ((A,\F)_1,\E_\Delta)_d = (\alpha \, f) (\alpha \, e) (f \, e)^{d-1} = 
(\alpha_\ux^2, (f \, e)^{d-1} \, f_\ux \, e_\ux)_2 = \\ 
& (A,(\F,\E_\Delta)_{d-1})_2 =0. \end{aligned} \] 
Since $(\F,\E_\Delta)_{d-1}$ is apolar to every order $2$ form, it must be identically zero. 

\subsection{} 
In fact we have an identity $(\F,\E_I)_{d-1}=0$ for any invariant. 
This can be informally explained as follows: $I$ is left unchanged by all 
$g \in SL_2$, hence it is annihilated by  the Lie algebra $\mathfrak{sl}_2$. Now 
observe that $\mathfrak{sl}_2$ (as the adjoint $SL_2$-representation) is 
isomorphic to $S_2$. The standard generators 
$\left( \begin{array}{cc} 0 & 1 \\ 0 & 0 \end{array} \right), 
\left( \begin{array}{cc} 0 & 0 \\ 1 & 0 \end{array} \right), 
\left( \begin{array}{rr} 1 & 0 \\ 0 & -1 \end{array} \right)$ respectively give the equations
(cf.~\cite[Theorem 4.5.2]{Sturmfels}) 
\[ \sum\limits_{i=0}^d \, (d-i) \, a_{i+1} \, \frac{\partial I}{\partial a_i}  = 
\sum\limits_{i=0}^d i \, a_{i-1} \, \frac{\partial I }{\partial a_i} = 
\sum\limits_{i=0}^d \, (d-2i)  \, a_i \, \frac{\partial I }{\partial a_i} = 0. \] 

\subsection{} \label{section.diffeq2.J}
Similarly we have a $(d-3)$-dimensional family of differential 
equations for $\Delta$ coming 
from the module $S_{d-4}$. We will express it in a form involving only the quadratic 
covariants of $\F$. As before, 
\[ ([(A,\Be_\F)_{d-4}]_{\uy = \ux},\E_\Delta)_d = 0 \quad \text{for all $A \in S_{d-4}$.} \] 
Let $A = {\alpha_\ux}^{d-4}, \, \Be_F = b_\ux^{d-2} \, {b'_\uy}^{d-2}$ where $b,b'$ are 
equivalent letters. Then 
\[ \begin{aligned} 
{} & ([(A,\Be_\F)_{d-4}]_{\uy = \ux},\E_\Delta)_d = ((\alpha \, b)^{d-4} \, 
{b_\ux}^2 \, {b_\ux'}^{d-2}, e_\ux^d \, )_d = \\
& (\alpha \, b)^{d-4} \, (b \, e)^2 (b' \, e)^{d-2} = 
(A,b_\ux^{d-4} \, (b \, e)^2 (b' \, e)^{d-2})_{d-4} =0, 
\end{aligned} \] 
hence 
\begin{equation} b_\ux^{d-4} \, (b \, e)^2 \, (b' \, e)^{d-2} =0. \label{bb'} 
\end{equation} 
Let $\ux \, \partial_\uy = 
x_1 \, \frac{\partial}{\partial y_1} + x_2 \, \frac{\partial}{\partial y_2}$, usually 
called the polarization operator. Then
$(\ux \, \partial_\uy)^2 \circ \Be_\F = 
(d-2) \, (d-1) \, b_\ux^{d-2} \, {b'_\ux}^2 \, {b'_\uy}^{d-4}$, 
hence identity (\ref{bb'}) is the same as 
\begin{equation} 
(\E_\Delta,(\ux \, \partial_\uy)^2 \circ \Be_\F)_d =0.
\label{quad-syzygy.2} \end{equation}
Let us write $(\F,\F)_{2r} = {\tau^{(2r)}_{\ux}}^{2d-4r}$ for the even quadratic 
covariants. 
We have a Gordan series (see~\cite[p.~55]{GrYo}) 
\[ \F(\ux) \, \F(\uy) = 
\sum\limits_{r=0}^{[\frac{d}{2}]} \, 
c_r \, (\ux \, \uy)^{2r} \, {\tau^{(2r)}_\ux}^{d-2r} \, {\tau^{(2r)}_\uy}^{d-2r}, \] 
where 
$c_r = \frac{\binom{d}{2r}^2}{\binom{2d-2r+1}{2r}}$. 
Apply the operator 
\[ \Omega = \frac{\partial^2}{\partial x_1 \, \partial y_2} - 
\frac{\partial^2}{\partial x_2 \, \partial y_1}, \] 
and divide by $(\ux \, \uy)$, then we get an expansion 
\[ \Be_\F = \frac{\Omega \circ \F(\ux) \, \F(\uy)}{(\ux \, \uy)} = 
\sum\limits_{r=1}^{[\frac{d}{2}]} \, c_r \, (2r) \, (2d-2r+1) \, 
(\ux \, \uy)^{2r-2} \, 
{\tau^{(2r)}_\ux}^{d-2r} \, {\tau^{(2r)}_\uy}^{d-2r}. \] 
Apply $(\ux \, \partial_\uy)^2$ to each term, which amounts to replacing 
the expression 
$(\ux \, \uy)^{2r-2} \, {\tau^{(2r)}_\ux}^{d-2r} \, {\tau^{(2r)}_\uy}^{d-2r}$ 
with 
\[ (d-2r)(d-2r-1) \, (\ux \, \uy)^{2r-2} \, 
{\tau^{(2r)}_\ux}^{d-2r+2} \, {\tau^{(2r)}_\uy}^{d-2r-2}. \] 
Now apply $(\E_\Delta,-)_d$ to each term, then 
\[ \begin{aligned} 
{} & (\epsilon_\ux^d, (\ux \, \uy)^{2r-2} \, 
{\tau^{(2r)}_\ux}^{d-2r+2} \, {\tau^{(2r)}_\uy}^{d-2r-2})_d = 
\epsilon_\uy^{2r-2} \, (\epsilon \, \tau)^{d-2r+2} \, {\tau^{(2r)}_\uy}^{d-2r-2} \\ 
& [(\epsilon_\ux^d, {\tau^{(2r)}_\ux}^{2d-4r})_{d-2r+2}]_{\ux = \uy}. \end{aligned} \] 
Hence finally we deduce the identity 
\begin{equation} \sum\limits_{r=1}^{[\frac{d-2}{2}]} \, 
\xi_r \, (\E_\Delta,(\F,\F)_{2r})_{d-2r+2} = 0, 
\label{diffeq.deg2} \end{equation}
where 
\[ \xi_r = \frac{(2d-4r+1)!}{(2r-1)!(d-2r-2)!(d-2r)!(2d-2r)!}. \] 

\subsection{The degree of the Jacobian scheme} 
Let $Z=\Proj (R/J)$. It is the scheme-theoretic degeneracy locus where the 
morphism 
\[ S_d \otimes \O_{\P^d} \lra S_2 \otimes \O_{\P^d}(1) \oplus 
S_{d-4} \otimes \O_{\P^d}(2) \] 
has rank $\le d-1$. Hence, by the Porteous formula (see~\cite[Ch.~II, \S 4]{ACGH})
the degree of $Z$ is given by the coefficient of $h^2$ in 
$(1+h)^{-3} (1+2h)^{3-d}$, which is $2 \, d \, (d-2)$. By Hilbert's formula 
in~\S\ref{deg.Xlambda}, 
\[ \deg X_{\tau} = 3 \, (d-2), \quad \deg X_{\delta} = 2 \, (d-2) \, (d-3). \] 
\begin{Proposition} \sl 
The scheme $Z$ has multiplicities $2$ and $1$ along $X_{\tau}$ and $X_{\delta}$ 
respectively. 
\label{prop.multiplicity} \end{Proposition} 
This means, for instance, that if $\eta_\tau$ is the scheme-theoretic 
generic point of $X_\tau$, then the ring $\O_{Z,\eta_\tau}$ is of length $2$. 

\smallskip 

\demo 
If the multiplicities are $a,b$, then 
$\deg Z = a \, \deg X_\tau + b \, \deg X_\delta$, i.e., 
\[ 2  d \, (d-2) = 3  a \, (d-2) + 2  b \, (d-2) \, (d-3). \] 
We have obvious constraints $a,b \ge 1$, and then it is straightforward 
to check that $(a,b)=(2,1)$ is the only possible solution. \qed 

\section{Examples} \label{section.examples} 
In this section we will describe $J$ and its primary decomposition for $d \le 5$. 
In each case the minimal system of generators for the ring of covariants was 
calculated in the nineteenth century (see~\cite[Ch.~V, VII]{GrYo}). 
If $C$ is a covariant of $d$-ics, then $\ga(C) \subseteq R$ will 
denote the graded ideal generated by the coefficients of $C$. 
\subsection{Cubics and quadratics} \label{d_2_or_3} 
So far we had assumed $d \ge 4$. The case $d=3$ is a little 
exceptional, but rather easy. In this case $Z$ is a non-reduced scheme 
of degree $6$ supported on the twisted cubic curve $X_{(3)}$. The minimal system 
of cubics consists of $\F,\He,\T$, and $\Delta = (\T,\F)_3$, i.e., every covariant is 
a polynomial function in these. It is immediate that $\E_\Delta= \T$, and 
Theorem~\ref{theorem.res.J} is true as stated with the convention that $S_{-1} =0$. 
Thus we have a resolution 
\[ 0 \la R/J \la R \la R(-3) \otimes S_3 \la R(-4) \otimes S_2 \la 0. \] 
The ideal of $X_{(3)}$ is $\ga(\He)$ (cf.~\cite[Exercise 11.32]{FH}), hence we have an 
equality $\ga(\He) = \sqrt{\ga(\T)}$. 

For $d=2$, we have $\Delta = (\F,\F)_2$ and $\E_\Delta = \F$, i.e., $J = (a_0,a_1,a_2)$ 
is the irrelevant maximal ideal. 

\subsection{Quartics} \label{d_4} 
Define $i = (\F,\F)_4, \, j = (\F,\He)_4$, 
which are invariants of degrees $2,3$. The minimal system for $d=4$ consists of 
$\F,\He,\T,i$ and $j$. Let $\gP_\tau,\gP_\delta \subseteq R$ denote the ideals 
of $X_{(3,1)}$ and $X_{(2,2)}$ respectively. 
\begin{Proposition} \sl 
\begin{enumerate} 
\item[(a1)]
We have identities 
\[ \Delta_\F = i^3 - 6 \, j^2, \quad 
\E_\Delta = i^2 \, F - 6 \, j \, \He. \] 
\item[(a2)]
$\gP_\tau$ is the complete intersection ideal $(i,j)$, and 
$\gP_\delta =\ga(\T)$. 
\item[(a3)] 
We have a primary decomposition 
$J = (i^2,j) \cap \gP_\delta$. 
\end{enumerate} 
\end{Proposition} 
\demo 
Since $\Delta$ is of degree $6$, it must be a linear combination of $i^3$ and 
$j^2$, say $c_1 \, i^3 + c_2 \, j^2$. Specialise to $F = x_1^2 \, x_2 \, (x_1 + x_2)$, 
when $\Delta_F$ must vanish. Computing directly, we get the equation 
$\frac{c_1}{216} + \frac{c_2}{1296} =0$, hence $c_1:c_2 = 1:-6$, i.e., 
we may take $\Delta = i^3 - 6 \, j^2$. Differentiating this identity, we get 
\[ \E_\Delta =  \frac{1}{6} \, (3 \, i^2 \times 2 \, \E_i - 
12 \, j \times 3 \, \E_j). \] 
But $\E_j = \He$ and $\E_i = \F$, hence it equals $i^2 \, \F - 6 \, j \, \He$. 
This proves (a1). 

Since $X_{(3,1)}$ is exactly the locus of nullforms, it is characterized by the 
vanishing of all invariants, i.e., $i=j=0$ at $F  \iff F \in X_{(3,1)}$. 
Since the ideal $(i,j)$ has no embedded primes, it must be $\gP_\tau$-primary. But since 
it also has degree $6 \, (= \deg \gP_\tau)$, we get $(i,j) = \gP_\tau$. 

In~\cite[Theorem 1.4]{AC1} it is proved that the ideal of every CR-locus of the 
type  $X_{(a,a)}$ is generated in degree $3$. It follows from the set-up described 
there that the degree $3$ piece $(\gP_\delta)_3$ is the kernel of the surjective 
morphism 
\[ S_3(S_4) \lra S_3(S_2 \otimes S_2) \lra S_3(S_2) \otimes S_3(S_2) \lra 
S_6 \otimes S_6 \lra S_2(S_6). \] 
We have plethysm decompositions 
\[ S_3(S_4) = S_{12} \oplus S_8 \oplus S_6 \oplus S_4 \oplus S_0, \quad 
 S_2(S_6) = S_{12} \oplus S_8 \oplus S_4 \oplus S_0, 
\] 
hence $(\gP_\delta)_3 \simeq S_6$. This subrepresentation must correspond to $\T$, 
since up to scalar it is the only covariant of degree-order $(3,6)$. 
This implies that $\gP_\delta = \ga(\T)$. 

To prove (a3), let $J = \gq_\tau \cap \gq_\delta$ be the (necessarily unique) primary 
decomposition, such that $\gq_\star$ is $\gP_\star$-primary.  
(See~\cite[Ch.~4]{AM} for generalities on primary decomposition.) 
Since $J$ has multiplicity one along $X_{(2,2)}$, we have 
$\gq_\delta = \gP_\delta$. Note that $(i^2,j)$ is $\gP_\tau$-primary 
(since it is perfect and its radical is $\gP_\tau$), moreover the 
expression for $\E_\Delta$ in (a1) shows that $J \subseteq (i^2,j)$. This implies 
that $\gq_\tau \subseteq (i^2,j)$, and it only remains 
to show the opposite inclusion. Let $z$ be any of the coefficients of $\T$, then 
\[ (J:z) = (\gq_\tau:z) \cap (\gP_\delta:z). \] 
Now $z \notin \gP_\tau$, hence $(\gq_\tau:z) = \gq_\tau$. Since $(\gP_\delta:z)=R$, we 
have $(J:z) = \gq_\tau$. From (a1), 
\[ (\E_\Delta,\He)_1 = (i^2 \, \F - 6 \, j \, \He,\He)_1 = 
i^2 \, (\F, \He)_1 - 6 \, j \, (\He,\He)_1 = i^2 \, \T,  \] 
and similarly $(\E_\Delta: \F)_1 = 6 \, j \, \T$.   
It follows that $i^2 \,z, j \, z \in J$, implying $i^2,j \in \gq_\tau$. 
This completes the proof of the proposition. \qed 

\medskip 

The identity~(\ref{diffeq.deg2}) of~\S\ref{section.diffeq2.J} 
reduces to $(\E_\Delta,\He)_4=0$, which gives the differential equation 
\[ \begin{aligned} 
{} & (2 \, a_0 \, a_2 - 2 \, a_1^2) \, \frac{\partial \Delta}{\partial a_0} + 
(a_0 \, a_3 - a_1 \, a_2) \, \frac{\partial \Delta}{\partial a_1} + 
(\frac{2}{3} \, a_1\, a_3 - a_2^2 + \frac{1}{3} \, a_0 \, a_4) 
\, \frac{\partial \Delta}{\partial a_2} + \\ 
& (a_1 \, a_4 - a_2 \, a_3) \, \frac{\partial \Delta}{\partial a_3} + 
(2 \, a_2 \, a_4 - 2 \, a_3^2) \, \frac{\partial \Delta}{\partial a_4}  =0. 
\end{aligned} \]

\subsection{Quintics} The invariant theory of the binary $d$-ic rapidly becomes more 
complicated with increasing $d$, in particular it is progressively harder to 
calculate $J$ precisely. In this section we will complete the calculation for 
$d=5$, making heavy use of machine computations in {\sc Maple} and Macaulay-2. 
The minimal system is given on~\cite[p.~131]{GrYo}. (Since it has 
$23$ members, it will not be reproduced here.) 
For quintics, the number of linearly independent covariants of degree-order $(m,q)$ is the 
number of copies of $S_q$ in the plethysm $S_m(S_5)$. We wrote 
our own set of {\sc Maple} procedures based on the Cayley-Sylvester formula 
(see ~\cite[Corollary 4.2.8]{Sturmfels}) to decompose it into irreducible summands. 

In addition to $\He$ and $\T$, we have covariants $i = (\F,\F)_4, A = (i,i)_2$ 
of degree-orders $(2,2),(4,0)$ respectively. Define 
\[ \begin{array}{lc} 
& \text{degree-order} \\ 
C_1 = 15 \, (i,\He)_2 + 2 \, i^2  & (4,4)\\
C_2 = 770 \, (i,\F \, \He)_2 - 675 \, (i,(\F,\He)_1)_1 + 
198 \, i^2 \, \F & (5,9) \\
D_1 = - 21  \, (C_1,\F^2)_4 + 55 \, (C_1,\He)_2 + 14 \, C_1 \, i  & (6,6) \\
D_2 = 5 \, (C_1,\He)_4 + 4 \, (C_1,i)_2 & (6,2) 
\end{array} \] 

\begin{Proposition} \sl 
\begin{enumerate}
\item[(b1)]
We have identities
\[ \begin{aligned} 
\Delta = & \, 59 \, A^2 + 320\,  (i^3,\He)_6, \\
\E_\Delta = & \, \frac{25}{3} \, A \, (i,\F)_1 +
\frac{3400}{21} \, i \, (i^2,\F)_3  \, - 240 \, (i^2,(\F,\He)_1)_4. 
\end{aligned} \] 
\item[(b2)]
If $\gP_\tau,\gP_\delta$ denote the ideals of $X_\tau, X_\delta$ respectively, 
then 
\[ \gP_\tau = \ga(C_1,A), \quad 
\gP_\delta = \ga(C_2). \] 
\item[(b3)]
We have a primary decomposition 
\[ J = \gq_\tau \cap \gP_\delta,  \] 
where $\gq_\tau = \ga(D_1,D_2)$ is $\gP_\tau$-primary. 
\end{enumerate} 
\end{Proposition}
\demo 
The minimal system shows that there are only two independent invariants 
in degree $8$, namely $A^2$ and $(i^3,\He)_6$. Hence $\Delta = 
c_1 \, A^2 + c_2 \, (i^3,\He)_6$ for some $c_i$. Specialise to
$F = x_1^2 \,  x_2 \, (x_1 + x_2) \, (x_1 - x_2)$ (when $\Delta$ must 
vanish), then we get $320 \, c_1 - 59 \, c_2 =0$. Similarly 
$A \, (i, \F)_1, i \, (i^2, \F)_3, (i^2,(\F,\He)_1)_4$ form a basis 
of covariants of degree-order $(7,5)$, hence $\E_\Delta$ must be 
their linear combination. We can find the coefficients by specialisation as before, and 
this establishes the formulae in (b1). 

First we determine the generators of $\gP_\tau$ using the recipe of 
\cite[\S 3.1]{ego3}. Write 
\[ \sum\limits_{i=0}^5 \, \binom{5}{i} \, a_i \, x_1^{5-i} \, x_2^i = 
(b_1 \, x_1 +  b_2 \, x_2)^3 \, 
(c_0 \, x_1^2 + 2 \, c_1 \, x_1 \, x_2 + c_2 \, x_2^2) \] 
(where $a,b,c$ are indeterminates), and equate the coefficients. This defines 
a ring morphism 
\[ \Complex [a_0,\dots,a_5] \lra \Complex[b_1,b_2,c_0,c_1,c_2],  \] 
whose kernel is $\gP_\tau$. A computation (done in Macaulay-2) shows 
that all the ideal generators are in degree $4$, and $\dim \,  (\gP_\tau)_4 = 6$. 
Now $A$ (being an invariant) must vanish on $X_\tau$, hence 
$(\gP_\tau)_4$ has $S_0$ as a summand. The module $S_4(S_5)$ contains no 
copies of $S_i$ for $0 < i < 4$, and $2$ copies of $S_4$. Hence 
$(\gP_\tau)_4$ must be isomorphic to $S_0 \oplus S_4$ as an 
$SL_2$-representation. The order $4$ piece (to be called $C_1$) must be a linear 
combination of $(i,\He)_2$ and $i^2$, because the latter form a basis in degree-order $(4,4)$. 
Then we determine the actual coefficients as before 
by specialising $F$ to $x_1^3 \, x_2 \, (x_1+x_2)$. 

A similar computation shows that $\gP_\delta$ is generated by a 
$10$-dimensional vector subspace of $R_5$. Notice that 
$X_\delta \supseteq X_{(4,1)}$, and by~\cite{Ei2}, the ideal of 
$X_{(4,1)}$ equals $\ga(i)$. Thus we have an inclusion $\gP_\delta \subseteq \ga(i)$; 
this implies that each degree $5$ covariant vanishing on $X_\delta$ must be 
a linear combination of terms of the form $(i, \Phi)_k$ for some 
degree $3$ covariant $\Phi$. (This follows because the vector space $(\ga(i))_5$ is spanned by 
such terms.) Clearly $0 \le k \le 2$. Now  
$S_3(S_5) \simeq S_{15} \oplus S_{11} \oplus S_9 \oplus S_7 \oplus S_5 \oplus 
S_3$, corresponding to the cases 
\[ \Phi = \F^3, \, \F \, \He, \, (\F,\He)_1, \, i \, \F, \, (i,\F)_1, \, (i,\F)_2. \] 
This allows us to write down all the possibilities for $(i,\Phi)_k$. 
An exhaustive search shows that $C_2$ is the only linear combination which vanishes 
on $F = x_1^2 \, x_2^2 \, (x_1 + x_2)$. This proves (b2). 

The $\gP_\delta$-primary component of $J$ is $\gP_\delta$ 
itself. Let $w$ denote the coefficient of $x_1^9$ in $C_2$, 
then $\gq_\tau$ (the $\gP_\tau$-primary component) equals the colon ideal 
$(J:w)$. We calculated the latter in Macaulay-2, and found it to have $10$ generators in 
degree $6$, and $12$ first syzygies in degree $7$. Hence we have a resolution 
\[ 0 \la R/\gq_\tau \la R \la R(-6) \otimes M_{10} 
\la R(-7) \otimes M_{12} \la \dots \] 
where $M_r$ denotes an $r$-dimensional $SL_2$-representation. Now 
\[ S_6 \, (S_5) = S_2^{\, \oplus 2} \oplus S_4 \oplus S_6^{\, \oplus 4} 
\oplus S_8^{\, \oplus 2} 
\oplus \text{summands $S_i$ with $i \ge 10$}, \] 
hence the dimension count forces $M_{10} \simeq S_6 \oplus S_2$. 
Let $D_1,D_2$ denote the corresponding covariants of orders $6$ and $2$. 
Since $\gq_\tau \subseteq \gP_\tau$, each $D_i$ can be written as a sum of terms 
of the form $(C_1,\Psi)_k, A \, \Psi'$, where $\Psi, \Psi'$ are of degree $2$. 
Thus we may write 
\[ \begin{aligned} D_1 & = \alpha_1  \, (C_1,\F^2)_4 + \alpha_2 \, 
(C_1,\He)_2 +\alpha_3 \, C_1 \, i, \\ 
D_2 & = \beta_1 \, (C_1,\He)_4 + \beta_2 \, (C_1,i)_2, \end{aligned} \] 
for some $\alpha_i, \beta_j \in {\mathbf Q}$. (The terms $A \, \He$ and $A \, i$ 
are not needed, because a calculation shows that they are respectively equal to 
\[ \frac{3}{25}  \, (C_1, \F^2)_4 - \frac{1}{25} \, (C_1,\He)_2 + 
\frac{162}{875} \, C_1 \, i, \quad  
\frac{18}{25} \, (C_1,\He)_4 -\frac{48}{125} \, (C_1,i)_2. ) \] 
Since $J \subseteq \gq_\tau$, we must have 
\[ \E_\Delta = \gamma_1 \, (D_1,F)_3 + \gamma_2 \, (D_2,F)_1 \] 
for some $\gamma_i \in {\mathbf Q}$. When rewritten in terms of 
the basis elements $A \, (i,\F)_1, i \, (i^2,\F)_3, (i^2,(\F,\He)_1)_4$ 
for covariants of degree-order $(7,5)$, this becomes an
inhomogeneous system of three linear equations. It turns out that there is a two-dimensional 
family of solutions, and the general solution can be written as 
\[ \begin{aligned} 
{} & (\gamma_1 \, \alpha_1,\gamma_1 \, \alpha_2,\gamma_1 \, \alpha_3,
\gamma_2 \, \beta_1,\gamma_2 \, \beta_2) = \\ 
& (-\frac{3}{5} - s + \frac{5}{4} \, t, \, 5 - \frac{5}{3} \, s - \frac{25}{6} \, t, 
- \frac{2}{7} - \frac{8}{7} \,  s  + \frac{75}{28} \,  t, \, s,t).  
\end{aligned} \] 
In order to determine $s,t$, we need to look at the first syzygies of $\gq_\tau$. Since they 
are all linear, $M_{12}$ must be a submodule of 
\[ M_{10} \otimes S_5 \simeq (S_6 \oplus S_2) \otimes S_5 \simeq 
S_{11} \oplus S_9 \oplus S_7^{\oplus 2} \oplus S_5^{\oplus 2} \oplus S_3^{\oplus 2} \oplus S_1.  \] 
By a dimension count, there are only four possible choices for $M_{12}$, it can only be 
$S_{11}, S_5^{\oplus 2}, S_5 \oplus S_3 \oplus S_1$ or $S_7 \oplus S_3$. It cannot be 
$S_{11}$ since the corresponding covariant is divisible by $\F$, and cancelling the 
latter would imply the absurdity that there is a first syzygy in degree $6$. 
If $S_5 \subseteq M_{12}$ 
(i.e., if there were a syzygy in order $5$), then there would be a nontrivial 
identity of the form $\eta_1 \, (D_1,F)_3 + \eta_2 \, (D_2,F)_1 = 0$. 
A calculation shows that there is none, this rules out all but the last choice. 
Thus $S_7 \subseteq M_{12}$, i.e., we have an identity of the form 
\[ \eta_1 \, (D_1,\F)_2 + \eta_2 \, D_2 \, \F = 0. \] 
Indeed, it turns out that 
$(s,t) = (\frac{24}{35}, \frac{96}{175}), \eta_1/\eta_2 = 4$ 
is the unique nontrivial solution. Finally we choose 
$\gamma_1 = \frac{1}{35}, \gamma_2 = \frac{24}{175}$, so that $D_1,D_2$ acquire 
integer coefficients. The proposition is proved. \qed 

\medskip 

It would be of interest to have a general result describing the primary decomposition of 
$J$ for all $d$, but this appears inaccessible. 

\subsection{} 
Not every invariant of binary forms has a perfect Jacobian ideal. E.g., let $d=4$ 
(with notation as in~\S \ref{d_4}). Let us show that $\gb = \ga(\E_j)$ 
(the Jacobian ideal of $j$) is not perfect. Since $\E_j$ is a covariant of 
degree-order $(2,4)$, it must coincide with $\He$ up to a scalar. The zero locus of 
$\gb = \ga(\He)$ is the rational normal quartic curve, 
hence $\dim \, (R/\gb) = 2$. However 
we have an identical relation $(\He,\F)_2 = \frac{1}{6} \, i \, \F$ (see~\cite[p.~92]{GrYo}),  
which implies that $i \, (a_0,\dots,a_4) \subseteq \gb$. Consequently 
$\gb$ is not a saturated ideal, and $\text{depth}\, (R/\gb) =0$.

\section{The binary resultant} \label{section.J_Res}
We begin with a recapitulation of the Cayley method of calculating 
the binary resultant (see \cite[Ch.~2]{GKZ}). The reader may also 
consult~\cite{Dickenstein-Andrea} for variations on this theme. Let 
\[ \F = \sum\limits_{i=0}^d \, \binom{d}{i} \, a_i \, x_1^{d-i} \, x_2^i, \quad 
\G = \sum\limits_{j=0}^e \, \binom{e}{j} \, b_j \, x_1^{e-j} \, x_2^j, \] 
denote generic binary forms of orders $d,e$. 
Define the product space 
$Y = \P S_d \times \P S_e \times \P S_1$ with projection maps $\mu_1,\mu_2,\pi$ onto 
the respective factors. Consider the subvariety 
\[ \tGamma = \{(F,G,l) \in Y: \text{$l$ divides $F,G$}\} \subseteq Y.  \] 
Let $f = \mu_1 \times \mu_2$, then 
$\Gamma = f(\tGamma) \subseteq \P^d \times \P^e$ is the resultant hypersurface. 

For any integers $m,n, p$, let $\O_Y(m,n,p)$ denote the line bundle 
\[ \mu_1^* \, \O_{\P^d}(m) \otimes \mu_2^* \, \O_{\P^e}(n) \otimes \pi^* \, \O_{\P^1}(p), \] 
with similar notation on $\P^d \times \P^e$. 
There is a tautological global section in $H^0(Y,\O_Y(1,0,d)) = S_d \otimes S_d$ 
corresponding to the trace element $\F$, and similarly for $\G$. 
Both of these sections simultaneously vanish at $(F,G,l)$ iff 
$(F,l^d)_d = (G,l^e)_e =0$, i.e., iff $l$ divides $F,G$. 
In fact we have a Koszul resolution 
\[ \begin{aligned} 
0 \ra \O_Y(-1,-1,-(d+e)) & \ra \O_Y(-1,0,-d) \oplus \O_Y(0,-1,-e) \\
& \ra \O_Y \ra \O_\tGamma \ra 0. 
\end{aligned} \] 
Now tensor with $\O_Y(0,0,d)$, and write this complex as 
\begin{equation} 
0 \ra \sC^{-2} \ra \sC^{-1} \ra \sC^0 \ra \O_{\tGamma}(0,0,d) \ra 0. 
\label{complex.C} \end{equation} 
We have a second quadrant spectral sequence 
\begin{equation} \begin{aligned} 
{} & E_1^{p,q} = R^q f_* \, \sC^p, 
\qquad d_r^{\, p,q}: E_r^{\, p,q} \ra E_r^{\, p+r,q-r+1}, \\
& E_\infty^{p+q} \Rightarrow R^{p+q} f_* \, \O_{\tGamma}(0,0,d) 
\end{aligned} \label{spectralseq} \end{equation} 
in the range $p=0,-1,-2$ and $q=0,1$. 

\subsection{} Now assume $d \ge e -1$, and $e \ge 2$. The only nonzero $E_1$ terms are 
\[ \begin{array}{ll} 
E_1^{-2,1} = \O(-1,-1) \otimes S_{e-2}, & E_1^{0,0} = \O \otimes S_d, \\ 
E_1^{-1,0} = \O(-1,0) \oplus \O(0,-1) \otimes S_{d-e}. 
\end{array} \] 
(Throughout $\O$ stands for $\O_{\P^d \times \P^e}$.) 
It is immediate that $R^i f_* \, \O_{\tGamma}(0,0,d) =0$ for $i > 1$, moreover we 
have exact sequences 
\[ \begin{aligned} 
{} & 0 \ra E_1^{-1,0} \ra E_1^{0,0} \ra E_2^{0,0} \ra 0, \\ 
& 0 \ra E_1^{-2,1} \stackrel{d_2^{-2,1}}{\ra} E_2^{0,0} 
\ra f_* \, \O_{\tGamma}(0,0,d) \ra 0. 
\end{aligned} \] 

\begin{Lemma} \sl 
The map $d_2^{-2,1}$ admits a unique $SL_2$-equivariant lifting (say $\theta$) 
to a map $E_1^{-2,1} \ra E_1^{0,0}$. 
\end{Lemma} 

\demo Indeed, the obstruction to this lift lies in the group 
\[ \text{Ext}^1(E_1^{-2,1},E_1^{-1,0}) = 
H^1(\O(0,1) \otimes S_{e-2}) \oplus H^1(\O(1,0) \otimes S_{e-2} \otimes 
S_{d-e})\] 
which is zero. Thus we have a surjection of $SL_2$-representations 
\begin{equation} \Hom(E_1^{-2,1}, E_1^{0,0}) \ra \Hom(E_1^{-2,1}, E_2^{0,0}). 
\label{surj.ss} \end{equation}
Since the construction of $d_2^{-2,1}$ is equivariant, it spans a copy of 
$S_0$ in the target of the map~(\ref{surj.ss}). By Schur's lemma 
it must come from an $S_0$ in the source, i.e., we have an equivariant lifting. 
If there were two such lifts, their difference would lie in 
\[ \begin{aligned} {} 
{} & \Hom (E_1^{-2,1}, E_1^{-1,0}) = 
H^0(\O(0,1)) \otimes S_{e-2} \oplus H^0(\O(1,0)) \otimes 
S_{e-2} \otimes S_{d-e} \\ 
& = \, [S_e \otimes S_{e-2}] \oplus [S_d \otimes S_{e-2} \otimes S_{d-e}].  
\end{aligned} \] 
However this is impossible; formula (\ref{Clebsch-Gordan}) from \S\ref{section.trans} 
shows that the last module does not contain {\sl any} copy of $S_0$. \qed 

\subsection{} 
Thus we get a map $E_1^{-2,1} \oplus E_1^{-1,0} \stackrel{\eta}{\lra} E_1^{0,0}$ of 
vector bundles of rank $d+1$ each, which can be seen as a map 
\begin{equation} 
S_{e-2} \oplus S_{d-e} \oplus S_0 \stackrel{\eta_{F,G}} \lra S_d 
\label{etaFG} \end{equation}
parametrised by points $(F,G) \in \P^d \times \P^e$. It fails to be bijective 
exactly over $\Gamma$. Now 
\[ \wedge^{d+1} \eta: \wedge^{e-1} \, E_1^{-2,1} \otimes 
\wedge^{d-e+2} \, E_1^{-1,1} \lra \wedge^{d+1} \, E_1^{0,0} \] 
is the map $\O(-e,-d) \lra \O$, i.e., $\R = \det \eta_{\F,\G}$ is an invariant of 
degree $(e,d)$ in the coefficients of $\F,\G$ respectively. Hence 
$\R$ must coincide with the resultant of $F,G$ (up to a scalar). 

The maps $S_0 \lra S_d, S_{d-e} \lra S_d$ are 
respectively $1 \ra F$, and $A \ra A \, G$ for $A \in S_{d-e}$. 
The map $\theta: S_{e-2} \lra S_d$ is given by the {\sl Morley form}  which we 
describe below. 
Symbolically write $F = f_\ux^d, \, G = g_\ux^e$. Define a joint covariant of 
$F,G$ by the expression 
\[ \Mor = \sum\limits_{i=1}^{e-1} \, (f \, g) \, 
f_\ux^{i-1} \, g_\ux^{e-i-1} \, f_\uy^{d-i} \, g_\uy^i.  \] 
It is of order $e-2$ and $d$ in $\ux,\uy$ respectively. 
\begin{Proposition} \label{prop.morley} 
\sl For $A = \alpha_\ux^{e-2} \in S_{e-2}$, the image $\theta(A)$ is given by 
\begin{equation} 
(-1)^{e-1} \, [(\Mor,A)_{e-2}]_{\uy= \ux} 
= - \, \sum\limits_{i=1}^{e-1} \, (f \, g) \,  
(\alpha \, f)^{i-1} \, (\alpha \, g)^{e-i-1} \, f_\ux^{d-i} \, g_\ux^i \, . 
\label{theta.A} \end{equation} 
\end{Proposition} 
The transvectant on the left hand side is with respect to $\ux$-variables, 
treating the $\uy$ as constants. 
The proof is postponed to \S\ref{proof.prop.morley}. 

\subsection{} \label{section.Res.resolution} 
Now the rest of the argument is very similar to the discriminant case. 
(At this point we leave the details to the reader.) 
That is to say, if $l \in S_1$ divides $F,G$, then each form in the image of the map 
\[ S_{e-2} \oplus S_{d-e} \lra S_d \] 
is divisible by $l$ (see~Lemma \ref{lemma.sigma.r} below), and the $\F$-evectant 
\[ \E_\R^{(\F)} = \sum\limits_{i=0}^d \, 
\frac{\partial \R}{\partial a_i} \, x_2^i \, (-x_1)^{d-i}, \] 
reduces to $\constant l^d$. In conclusion, we get the following result: 

\begin{Theorem} \sl 
The ideal 
\[ J_\F = (\frac{\partial \R}{\partial a_0}, \dots, \frac{\partial \R}{\partial a_d})
\subseteq Q = \Complex[a_0,\dots,a_d,b_0,\dots,b_e] \] 
is perfect of height $2$, with an equivariant bigraded minimal resolution 
\[ \begin{aligned} 
0 \la Q/J_\F \la & Q \la  Q(1-e,-d) \otimes S_d \la \\ 
& Q(1-e,-d-1) \otimes S_{d-e} \oplus Q(-e,-d-1) \otimes S_{e-2} \la 0. 
\end{aligned} \] 
\label{EF.perfect} \end{Theorem}

\subsection{} The syzygy modules $S_{d-e}$ and $S_{e-2}$ respectively 
correspond to the identities 
\[ (\G,\E_\R^{(\F)})_e =0, \quad 
   (\Mor,\E_\R^{(\F)}|_{\uy=\ux})_d^\uy =0. \] 
In the latter, we have changed $\E$ into a $\uy$-form of order $d$. The 
transvection is with respect to $\uy$-variables, leaving an $\ux$-form 
of order $e-2$. We will rewrite this identity non-symbolically, in a form which 
only involves the joint covariants $(\F,\G)_r$. First we expand 
each term of $\Mor$ into its Gordan series (see \cite[p.~55]{GrYo}), 
i.e., we write 
\begin{equation} f_\ux^{i-1} \, g_\ux^{e-i-1} \, f_\uy^{d-i} \, g_\uy^i = 
\sum\limits_{s=0}^{e-2} \, \alpha_s \, (\ux \, \uy)^s \, 
(\uy \partial_\ux)^{d-s} \circ 
[(f_\ux^{i-1} \, g_\ux^{e-i-1}, f_\ux^{d-i} \, g_\ux^i)_s],
\label{Gseries} \end{equation}
where 
\[ \alpha_s = \frac{\binom{d}{s} \, \binom{e-2}{s}}
{\binom{d+e-s-1}{s}\binom{d+e-2s-2}{d-s}}. \] 
Using the general formalism of~\cite[\S 3.2.5]{Glenn}, 
\[ (f_\ux^{i-1} \, g_\ux^{e-i-1}, f_\ux^{d-i} \, g_\ux^i)_s = 
\beta_{i,s} \, (f \, g)^{s+1}  \, f_\ux^{d-s-1} \, g_\ux^{e-s-1}, \] 
where 
\[ \beta_{i,s} = \frac{1}{\binom{d}{s} \binom{e-2}{s} s!} 
\sum\limits_{l=0}^s (-1)^l \, l! (s-l)! \binom{i-1}{s-l} 
\binom{e-i-1}{l} \binom{d-i}{l} \binom{i}{s-l}. \] 
Now $(f \, g)^{s+1} \, f_\ux^{d-s-1} \, g_\ux^{e-s-1} = 
(\F,\G)_{s+1}$, which we write symbolically as $\tau_\ux^{d+e-2s-2}$. 
Then 
\[ (\uy \partial_{\ux})^{d-s} \, \circ \tau_\ux^{d+e-2s-2} = 
\binom{d+e-2s-2}{d-s} \, \tau_\ux^{e-s-2} \, \tau_\uy^{d-s}. \] 
Writing $\E_\R^{(\F)}|_{\uy = \ux} = \epsilon_\uy^d$, 
\[ \begin{aligned} 
{} & ((\ux \, \uy)^s \, \tau_\ux^{e-s-2} \, \tau_\uy^{d-s}, 
\epsilon_\uy^d)_d^\uy =  
(-1)^s \, \epsilon_\ux^s \, \tau_\ux^{e-s-2} \, (\tau \, \epsilon )^{d-s} \\
= \; & (-1)^d \, (\E_\R^{(\F)},(\F,\G)_{s+1})_{d-s}. 
\end{aligned} \] 
Hence, by substituting into~(\ref{Gseries}) we get the required identity 
\begin{equation}
\sum\limits_{s=0}^{e-2} \, \omega_s \, ( \E_\R^{(\F)},(\F,\G)_{s+1})_{d-s} =0, 
\end{equation}
where 
$\omega_s = \binom{d+e-2s-2}{d-s} \, \alpha_s \, \sum\limits_{i=1}^{e-1} \beta_{i,s}$. 
\subsection{} 
If $e=1$, then Theorem~\ref{EF.perfect} is true as stated if we take $S_{-1} =0$. 
If $d -e < -1$, then the spectral sequence~(\ref{spectralseq}) has a 
nonzero term at $E_1^{-1,1}$. We still get a determinantal formula 
\[ \R = \det \, (S_{e-2} \oplus S_0 
\stackrel{\eta'_{\F,\G}}{\lra} S_d \oplus S_{e-d-2}), \] 
but $J$ may no longer be perfect. E.g., for $(d,e)=(2,4)$, a Macaulay-2 computation shows that 
$J$ is of height $2$, but $\text{proj-dim}_Q \, (Q/J_{\F})=3$. 

\subsection{} \label{proof.prop.morley} 
Now we take up the proof of Proposition~\ref{prop.morley}. 
For $i=1,2$, let $U_i = \{ l \in S_1: \frac{\partial \, l}{\partial x_i} \neq 0\} \subseteq \P^1$, and 
$\U_i = \pi^{-1}(U_i)$. 
We will calculate the differential $d_2^{-2,1}$ using a {\v C}ech resolution of the 
complex~(\ref{complex.C}) for the cover $\U_i$. 
Let us write $\Sz_k^j$ as an abbreviation for $f_*(\sC^j|_{\U_k})$, where 
$k$ may denote $1,2$, or $12$. (As usual $\U_{12} = \U_1 \cap \U_2$.) On 
$\P^d \times \P^e$ we have a double complex of locally free sheaves
\[ \diagram 
\Sz^{-2}_{12} \rto^{h_1} & \Sz^{-1}_{12} \rto & \Sz^{0}_{12} \\ 
\Sz^{-2}_1 \oplus \Sz^{-2}_2 \rto \uto & \Sz^{-1}_1 \oplus \Sz^{-1}_2 \rto^{h_3} \uto^{h_2} & 
\Sz^{0}_1 \oplus \Sz^{0}_2 \uto \enddiagram \] 
It will be convenient to see it as a diagram of morphisms of vector spaces 
parametrised by the pair $(F,G)$. 
Since expression~(\ref{theta.A}) is linear in $A$, it is enough to show the proposition 
for a monomial $A$. Let $A = x_1^r \, x_2^{e-2-r}$. 

The isomorphism $S_{e-2} \simeq S_{e-2}^*$ of \S\ref{section.trans} takes 
the form $A$ to 
$A' = (-1)^{e-2-r} \, \binom{e-2}{r} \, x_2^r \, x_1^{e-2-r}$, 
since $(A,A')_{e-2} =1$. 
This implies that the sequence of isomorphisms 
\[ S_{e-2} \simeq S_{e-2}^* \simeq H^0(\P^1,\O(e-2))^* \otimes H^1(\P^1,\O(-2)) 
\simeq H^1(\P^1,\O(-e)), \] 
takes $A$ to the {\v C}ech cocyle 
\[ \bA = \frac{1}{A'} \, \times \frac{1}{x_1 \, x_2} = 
\frac{(-1)^{e-2-r}}{\binom{e-2}{r} \, x_1^{e-1-r} \, x_2^{r+1}} \in H^0(U_{12},\O(-e)).  \] 
Recall that by the usual procedure for calculating the differentials in a spectral 
sequence (see~\cite[\S 14]{BottTu}), 
\[ d_2^{-2,1}(\bA) = h_3 \circ h_2^{-1} \circ h_1(\bA). \] 
(Throughout, the vector space morphisms over $(F,G)$ are also denoted by $h_i$.)  
\subsection{} \label{section.FGrGFr} 
By the construction of the Koszul complex, $h_1(\bA) = F \, \bA \oplus G \, \bA$. 
To take the pre-image by $h_2$, 
we need to rewrite each of the summands as a difference $e^{(1)}-e^{(2)}$, 
where the denominator of $e^{(i)}$ is a power of $x_i$ alone. Write 
\[ F = \frac{(d-e+1)!}{d!} \, 
[ \, (y_1 \frac{\partial}{\partial x_1} + y_2 \, \frac{\partial}{\partial x_2})^{e-1} \, 
F \, ]_{\uy = \ux}.  \] 
Expand and retain only those terms whose power in $x_1$ is at least $e-1-r$, i.e., 
let 
\[ {\hat F} = \frac{(d-e+1)!}{d!} \sum\limits_{q \ge e-1-r} 
\, \binom{e-1}{q} \, x_1^q \, x_2^{e-q-1} \, 
\frac{\partial^{e-1} F}{\partial x_1^q \, \partial x_2^{e-q-1}}. \] 
Multiplying by $\bA$, we get 
\begin{equation} \begin{aligned} 
{} & e^{(2)} = \tF_r \\ 
& = \frac{(-1)^{e-2-r} \, (d-e+1)!}{x_2^{r+1} \, \binom{e-2}{r} \, d!}
\sum\limits_{q=e-r-1}^{e-1} \binom{e-1}{q} \, 
x_1^{q-e+r+1} \, x_2^{e-q-1} \, 
\frac{\partial^{e-1} F}{\partial x_1^q \, \partial x_2^{e-q-1}}, 
\end{aligned} \label{tFr} \end{equation} 
and then $e^{(1)} = F - \tF_r$. Similarly, let 
\begin{equation} 
\tG_r = \frac{(-1)^{e-2-r}}{x_2^{r+1} \, \binom{e-2}{r} \, e!} \sum\limits_{q=e-r-1}^{e-1}
\binom{e-1}{q} \, x_1^{q-e+r+1} \, x_2^{e-q-1} \, 
\frac{\partial^{e-1} G}{\partial x_1^q \, \partial x_2^{e-q-1}} 
\label{tGr} \end{equation}
Now $u = (F-\tF_r,-\tF_r) \oplus (G-\tG_r,-\tG_r)$ is an element such that 
$h_2(u) = F \bA \oplus G \bA$. 

To calculate the image of $u$ by $h_3$, multiply the first summand by $G$, the second 
by $F$ and subtract. This gives 
\begin{equation} d_2^{-2,1}(\bA) = h_3(u) = F \, \tG_r - G \, \tF_r. 
\label{d2A} \end{equation}
(Note that we have used a hidden `term order' where $F$ comes before $G$. As long as 
we remain consistent, this should cause no harm.) 

It is not \emph{a priori} obvious that the result is invariant 
under a change of variables, since the {\v Cech} cover is clearly not so invariant. 
On the other hand, expression~(\ref{theta.A}) is entirely in terms of symbolic brackets, 
hence visibly invariant. Thus, to complete the proof, we have to establish the identity 
\begin{equation} 
F \, \tG_r - G \, \tF_r = (-1)^{e-1} \, [(\Mor,A)_{e-2}]_{\uy= \ux}. 
\label{morley.identity} \end{equation} 
This calculation is done in the appendix. 

The following lemma was needed in \S\ref{section.Res.resolution}. 
\begin{Lemma} \sl 
If $l$ divides $F$ and $G$, then it divides $\theta(A)$ for any 
$A$. 
\label{lemma.sigma.r} \end{Lemma} 
\demo We may assume that $l = x_1$. Since $\theta$ is linear, it suffices to 
give a proof for a monomial $A$. But then the claim follows 
because $x_1$ clearly divides the left hand side of~(\ref{morley.identity}). \qed

\section{The $\Phi_n$ are arithmetically Cohen-Macaulay} \label{section.acm}
Let $\Phi_n \subseteq \P S_d$ be as in the introduction, with ideal 
$I_n \subseteq R$. We will exhibit $\Phi_n$ as the degeneracy locus of a map of 
vector bundles and then deduce that $I_n$ is perfect ideal. Along the way 
we will construct a covariant $\A_n$ of binary $d$-ics such that 
$F \in \Phi_n \iff \A_n(F)=0$. 

\subsection{} Every $F \in S_d$ has a factorization 
\begin{equation} F = l_1^{ e_1} \dots l_n^{e_n} 
\label{F.factors} \end{equation} 
where the $l_i$ are pairwise nonproportional, and 
$e_1 \ge \dots \ge e_n > 0$. Let 
\[ g_F = \gcd \, (F_{x_1},F_{x_2}). \] 

\begin{Lemma} \sl 
With notation as above, $g_F = \prod\limits_i \, l_i^{ e_i-1}$. 
\label{lemma.gcd} \end{Lemma} 
\demo Evidently $g = \prod l_i^{ e_i-1}$ divides both the 
$F_{x_i}$, write $F_{x_1} = g \, A, F_{x_2} = g \, B$. Divide 
Euler's equation $d \, F = x_1 \, F_{x_1} + x_2 \, F_{x_2}$ 
by $g$, then 
$d \, \prod l_i = x_1 \, A + x_2 \, B$. If $A,B$ have a common linear 
factor, it must be one the $l_i$, say $l_1$. But 
\[ A = \sum\limits_i \, e_i \frac{\partial l_i}{\partial x_1} 
(\prod\limits_{j \neq i} l_j), \] 
so $l_1|A$ implies $\frac{\partial l_1}{\partial x_1}=0$. The same 
argument on $B$ leads to $\frac{\partial l_1}{\partial x_2}=0$, 
so $l_1=0$. This is absurd, hence $A,B$ can have no common factor, 
i.e.,~$g = g_F$. \qed 

\begin{Corollary} \sl 
Let $F \in S_d$. Then $F \in \Phi_n$ iff $\ord g_F \ge d-n$. \qed 
\end{Corollary} 

\subsection{} 
We have a map $S_d \otimes S_1 \lra S_{d-1}$ by formula~(\ref{Clebsch-Gordan}), 
we may see it as a morphism of vector bundles 
$\O_{\P^d}(-1) \otimes S_1 \lra S_{d-1}$.  Now consider the 
composite 
\[ \O_{\P^d}(-1) \otimes S_1 \otimes S_{n-1} \lra S_{d-1} \otimes S_{n-1} 
\stackrel{\text{mult}}{\lra} S_{d+n-2}, \] 
which we denote by $\alpha_n$. On the fibres over 
$[F] \in \P^d$, this can be thought of as a morphism 
\[ \begin{aligned} \alpha_{n,F}: \, 
S_1 & \otimes S_{n-1} \lra S_{d+n-2}, \\ 
l & \otimes G \lra (l,F)_1 \, G = \constant \, (l_{x_1}F_{x_2}-l_{x_2}F_{x_1}) \, G.
\end{aligned} 
\] 
Now $F_{x_1},F_{x_2}$ have a common factor of order $\ge d-n$, iff 
there are order $n-1$ forms $G_1,G_2$ such that 
$G_2 \, F_{x_1} + G_1 \, F_{x_2}=0$. This condition can be rewritten as 
$\alpha_{n,F}(x_1 \otimes G_1 - x_2 \otimes G_2)=0$. Hence $\alpha_{n,F}$ fails to be 
injective iff $F \in \Phi_n$. 

Let $\Psi_n$ denote the determinantal 
scheme $\{ \rank (\alpha_n) < 2n \}$ locally defined by the maximal 
minors of the matrix of $\alpha_{n,F}$. We have shown that 
$(\Psi_n)_\text{red} = \Phi_n$. 

\begin{Theorem} \sl 
The scheme $\Psi_n$ is reduced, hence $\Psi_n = \Phi_n$ as schemes. 
\end{Theorem}
\demo 
The standard codimension estimate for determinantal loci 
(see \cite[Ch.~2]{ACGH}) takes the form 
\[ \codim \Psi_n \le d+n-1-(2n-1) = d-n. \] 
Since equality holds, $\Psi_n$ is a Cohen-Macaulay scheme, in 
particular it has no embedded components. By the Thom-Porteous 
formula, 
\[ \deg \Psi_n = 
(-1)^{d-n} \times \text{coefficient of $h^{d-n}$ in} \; 
(1-h)^{d+n-1} = \binom{d+n-1}{d-n}. \] 
If we show that this coincides with $\deg \Phi_n$, then 
it will follow that $\Psi_n$ is reduced. 

Let $\lambda$ be a partition of $d$ with $n$ parts. 
A moment's reflection will show that $\deg X_\lambda$ as given by 
Hilbert's formula is the coefficient of the monomial 
$\prod\limits_{r=1}^d \, {z_r}^{r e_r}$ in the expression
\[ (z_1 + 2 \, z_2^2 + \dots + r \, z_r^r + \dots )^n. \] 
Now substitute the same letter $z$ for each $z_r$, then $\prod {z_r}^{r e_r} = z^d$. 
Hence the coefficient of $z^d$ in $(z + 2 \, z^2 + \dots + r \, z^r + \dots )^n$
equals 
\[ \sum\limits_{\text{$\lambda$ has $n$ parts}} \deg X_\lambda = 
\deg \Phi_n. \] 
But 
\[ (z + 2 \, z^2 + \dots + r \, z^r + \dots ) = \frac{z}{(1-z)^2},  \] 
hence this coefficient is the same as 
\[ \text{coefficient of $z^{d-n}$ in $(1-z)^{-2n}$}
= (-1)^{d-n} \, \binom{-2n}{d-n} = \binom{d+n-1}{d-n}. \] 
This completes the proof of the theorem. \qed

\medskip 

It follows that the Eagon-Northcott complex of the map 
\[ R(-1) \otimes S_1 \otimes S_{n-1} \lra R \otimes S_{d+n-2} \] 
gives a resolution of $R/I_n$ (see \cite[Ch.~2C]{BrunsVetter}). Its terms are: $\E^0 = R$, and 
\begin{equation} \E^p = \wedge^{2n-p-1} (S_{d+n-2}) \otimes 
S_{-(p+1)}(S_1 \otimes S_{n-1}) \otimes R(-2n+p+1), 
\label{ENcomplex} \end{equation} 
for $-(d-n) \le p \le -1$. 

\subsection{The covariants $\A_n$}
Consider the map 
\[ \wedge^{2n} \, \alpha_{n,\F}: \Complex \lra 
\wedge^{2n} \, S_{d+n-2}. \] 
Let $\A_n$ denote the image of $1$ via this map, which is a covariant of degree-order 
$(2n,2n(d-n-1))$ of binary $d$-ics. (It is well-defined only up to a multiplicative constant.) 
By construction, it is the Wronskian of the forms 
\begin{equation}
 \{ x_1^{n-j-1} \, x_2^j \, F_{x_i}: 0 \le j \le n-1, i =1,2 \}, 
\label{im-span} \end{equation}
i.e., it is the determinant of the following $2n \times 2n$ matrix: 
\begin{equation} (p,q) \lra 
\begin{cases} 
(x_1^{2n-q-1} \, x_2^q, \, x_1^{n-p-1} \, x_2^p \; \F_{x_1})_{2n-1} 
& \text{if $ 0 \le p \le n-1$,} \\ 
(x_1^{2n-q-1} \, x_2^q, \, x_1^{2n-p-1} \, x_2^{p-n} \; \F_{x_2})_{2n-1} 
& \text{if $ n \le p \le 2n-1$,} 
\end{cases} \label{formula.An} \end{equation}
and $0 \le q \le 2n-1$. It vanishes at $F$ iff the collection~(\ref{im-span}) is 
linearly dependent, hence 
\begin{Corollary} \sl 
$F \in \Phi_n \iff \A_n(F)=0$. 
\end{Corollary} 
Since $\A_{d-1}$ is an invariant of degree $2(d-1)$ it must coincide with the discriminant. 
Similarly $\A_1$ is (up to a scalar) the same as the Hessian. 
Thus the series $\{A_n\}$ can be thought of as an `interpolation' between the two. 

The following lemma will be used in the next section. 
\begin{Lemma} \label{lemma.gF} \sl 
Asssume $[F] \in \Phi_n \setminus \Phi_{n-1}$. Then 
\[ \A_{n-1}(F) = (g_F)^{2n-2}. \] 
\end{Lemma} 
\demo This is perhaps best proved using the relation between the 
Wronskian and ramification indices (see~\cite[pp.~37--43]{ACGH}). By hypothesis, 
$\alpha_{n-1}$ is of rank $2n-2$ at $[F]$, in fact 
\[ \im(\alpha_{n-1,F}) = \{ F_{x_1} \, G_1 + F_{x_2} \, G_2: G_i \in S_{n-2} \} 
= \{ g_F \, G: G \in S_{2n-3} \}. \] 
This can be seen as a linear series $\Sigma$ on $\P^1$ of degree 
$d+n-2$ and dimension $2n-3$. Write $F = \prod l_i^{e_i}$, then $\Sigma$ 
is only ramified at points $p_i \in \P^1$ corresponding to the $l_i$. 
Its ramification indices at $p_i$ are 
\[ e_i, e_i + 1, \dots, e_i + 2n-3. \] 
Hence the Wronskian of $\im(\alpha_{n-1,F})$ is 
$\prod {l_i}^{(2n-2) e_i} = (g_F)^{2n-2}$.  \qed 

\subsection{The codimension two case} 
Assume $n=d-2$. 
Then in the complex~(\ref{ENcomplex}) we have 
\[ \E^{-1} = \wedge^{2d-4} \, S_{2d-4} \otimes R(-2d+4) = 
S_{2d-4} \otimes R(-2d+4), \] 
i.e., $I_{d-2} = \ga(\A_{d-2})$. Now $J$ (the Jacobian ideal of $\Delta$) is contained in 
the ideal $I_{d-2}$, hence the image of the natural multiplication map 
\[ (I_{d-2})_{2d-4} \otimes R_1 \lra R_{2d-3} \] 
must contain the representation $(J)_{2d-3}$. Since the latter is spanned by the 
coefficients of $\E_\Delta$, we deduce the following: 
\begin{Corollary} \sl 
The covariants $(\A_{d-2},\F)_{d-2}$ and $\E_\Delta$ are equal up to a nonzero scalar. 
\end{Corollary} 

We end this section by constructing covariants which distinguish between 
the components $X_{\tau} = X_{(3,1^{d-3})}$ and $X_{\delta}=X_{(2^2,1^{d-4})}$. 
A result due to Hilbert~\cite{Hilbert1} says that a binary $d$-ic $F$ lies in 
$X_{(d)}$ iff $\He(F)=0$, and it lies in $X_{(d/2,d/2)}$ (assuming $d$ even) 
iff $\T(F)=0$. 

First assume that $F \in X_\tau \setminus X_\delta$. Then 
$g_F = l^2$ for some $l \in S_1$, and then $\A_{d-3} = l^{4d-12}$ by Lemma~\ref{lemma.gF}. 
If $F \in X_\delta \setminus X_\tau$, then 
$g_F = l_1 \, l_2$ for some nonproportional linear forms, and 
$\A_{d-3} = (l_1 \, l_2)^{2d-6}$. Hence we get the following proposition. 
\begin{Proposition} \sl 
Let $F$ be a binary $d$-ic. Then 
\[ \begin{aligned} 
F \in X_\tau  & \iff \A_{d-2}(F) = \He(\A_{d-3}(F)) =0, \\ 
F \in X_\delta & \iff \A_{d-2}(F) = \T(\A_{d-3}(F)) =0. 
\end{aligned} \] 
\end{Proposition} 

\medskip 

\begin{Remark} \rm 
Throughout this paper we have used $\Complex$ as our base field. Note however, that all the 
irreducible representations of $SL_2 \, {\mathbf Q}$ are defined over ${\mathbf Q}$, 
hence so are all the varieties and schemes defined above. Thus all of our results are 
valid over an arbitrary field of characteristic zero. 
\end{Remark} 

\newpage 

\section{Appendix: the Morley form} 
\vspace{-0.2cm} 
\centerline{(by A.~Abdesselam)} 
\vspace{0.2cm} 

\subsection{} We will now prove identity~(\ref{morley.identity}) 
from~\S\ref{section.FGrGFr}. 
At this point, a brief explanatory remark on the symbolic method should be 
helpful. 
We have $f_\ux = (f_1 \, x_1 + f_2 \, x_2), g_\ux = (g_1 \, x_1 + g_2 \, x_2)$ 
where $f_i,g_i$ are treated as indeterminates. Introduce the differential operators 
\[ \DD_F = \frac{1}{d!} \, F(\frac{\partial}{\partial f_1},\frac{\partial}{\partial f_2}), 
\quad \DD_G = \frac{1}{e!} \, G(\frac{\partial}{\partial g_1},\frac{\partial}{\partial g_2}). \] 
Then we have identities $F = \DD_F \, f_\ux^d, G = \DD_G\, g_\ux^d$. Moreover, 
each well-formed symbolic expression in $f,g$ can be evaluated by subjecting 
it to these operators; this is one way of providing a rigorous justification for 
the method. Thus the Morley form will be written as 
\[ \Mor(\ux,\uy) =  \DD_F \, \DD_G \, 
\sum\limits_{i=1}^{e-1} \, (f \, g) \, f_\ux^{i-1} \, g_\ux^{e-i-1} \, f_\uy^{d-i} \, g_\uy^i . 
\label{eq4} \] 
Now let 
\[ \theta_r = (-1)^{e-1} \, 
[ \, (\Mor, x_1^r \, x_2^{e-2-r})_{e-2} \, ]_{\uy:=\ux} \, ,  \] 
where the transvection is with respect to $\ux$. By definition, 
\[ \theta_r = 
\frac{(-1)^{e-1}}{(e-2)!^2}  \left. \left\{
(\frac{\partial^2}{\partial z_1 \, \partial x_2}
-\frac{\partial^2}{\partial z_2 \, \partial x_1})^{e-2} 
\,  \Mor(\uz,\uy) \, x_1^r \, x_2^{e-2-r}  \right\} \right|_{\uy:=\ux .} \] 
After a binomial expansion this simplifies to 
\begin{equation} 
\frac{(-1)^{e-1}}{(e-2)!} \, 
(-\frac{\partial}{\partial z_2})^r \, (\frac{\partial}{\partial z_1})^{e-2-r} \, 
\Mor(\uz,\ux). \label{expression.theta.r}
\end{equation}

\subsection{} Let us introduce a pair of variables 
$b = (b_1,b_2)$, which will serve as placeholders. Define the sum 
\begin{equation}  \Psi=(-1)^e \, \sum_{r=0}^{e-2} \, \binom{e-2}{r} \, 
b_1^r \, b_2^{e-2-r} \, (F \, \tG_r - G \, \tF_r), 
\label{Psi.defn} \end{equation}
so that 
\begin{equation} 
F \, \tG_r - G \, \tF_r = 
\frac{(-1)^e}{(e-2)!} \, \frac{\partial^{e-2} \, \Psi}
{\partial \, b_1^r \, \partial \, b_2^{e-2-r}}
\label{FGr-GFr}  \end{equation} 
Now 
\[ \begin{aligned} 
{} & G \, \frac{\partial^{e-1}F}{\partial x_1^q \, \partial x_2^{e-1-q}}  = 
(\DD_G \; g_\ux^e) \, 
[ (\frac{\partial^{e-1}}{\partial x_1^q \, \partial x_2^{e-1-q}}) \, \DD_F \; f_\ux^d \, ] \\
& =\DD_F \, \DD_G \, [ \, g_\ux^e  \, 
(\frac{\partial^{e-1}}{\partial x_1^q \, \partial x_2^{e-1-q}}) \, f_\ux^d \, ] \\ 
& = \frac{d!}{(d-e+1)!} \, \DD_F \, \DD_G \left[ f_1^q \, f_2^{e-1-q} \, 
f_\ux^{d-e+1} g_\ux^e \, \right], 
\end{aligned} \] 
and similarly 
\[ F \, \frac{\partial^{e-1}G}{\partial x_1^q \, \partial x_2^{e-1-q}}= 
e! \,  \DD_F \, \DD_G \, [ \, g_1^q \, g_2^{e-1-q} \, f_\ux^d \, g_\ux \, ].  \] 
Now substitute these expressions into equations~(\ref{tFr}) and (\ref{tGr}) 
from~\S\ref{section.FGrGFr}, and 
then substitute the latter into (\ref{Psi.defn}). Then we have 
$\Psi = \DD_F \, \DD_G \, \tPsi$, where
\[ \begin{aligned} 
{} & \tPsi = \sum\limits_{r=0}^{e-2} \; \big[ (-b_1)^r \, b_2^{e-2-r} \, \times  \\ 
& \sum\limits_{q=e-r-1}^{e-1} \binom{e-1}{q} \, x_1^{q-e+r+1} \, x_2^{e-q-r-2} 
\{ g_1^q \, g_2^{e-1-q} \, f_\ux^d \, g_\ux  - f_1^q \, f_2^{e-1-q} \, f_\ux^{d-e+1} \, g_\ux^{e} \}
\big]. \end{aligned} \] 
The double sum is over the range 
$0 \le r \le e-2, \, e-r-1 \le q \le e-1$, which is the same as 
$1 \le q \le e-1, \, e-q-1 \le r \le e-2$. 
Therefore, after changing the order of summation,
\[ \begin{aligned} \tPsi = 
& \sum\limits_{q=1}^{e-1} \, \big [\binom{e-1}{q} b_2^{e-2} x_1^{q-e+1} x_2^{e-q-2} 
\{ g_1^q \, g_2^{e-1-q} \, f_\ux^d \, g_\ux  - 
f_1^q \, f_2^{e-1-q} \, f_\ux^{d-e+1} \, g_\ux^e \} \times \\ 
& \sum\limits_{r=e-1-q}^{e-2} (-\frac{b_1 \, x_1}{b_2 \, x_2})^r \big], 
\end{aligned} \] 
which we abbreviate to 
\[ \sum\limits_{q=1}^{e-1} \, [ \, (\zA_1 - \zA_2) \times 
\sum\limits_{r=e-1-q}^{e-2} (-\frac{b_1 \, x_1}{b_2 \, x_2})^r \big]. \] 
The geometric series over $r$ is equal to 
\[ \begin{aligned} 
{} & \; \frac{(-b_1 \, x_1)^{e-1-q} \, (b_2 \, x_2)^q - 
(-b_1 \, x_1)^{e-1}}{(b_2 \, x_2)^{e-2} \, b_\ux} \\ 
= & \; (-b_1 \, x_1)^{e-1-q} \, (b_2 \, x_2)^{q-e+2} \, b_\ux^{-1} - 
(-b_1 \, x_1)^{e-1} \, (b_2 \, x_2)^{-e+2} \, b_\ux^{-1} \\
= & \; \zB_1 - \zB_2. 
\end{aligned} \] 
Hence, after expansion $\tPsi$ is a sum of four terms 
\[ \underbrace{\sum \zA_1 \, \zB_1}_{T_1}  + 
\underbrace{\sum - \zA_1 \, \zB_2}_{T_2} + \underbrace{\sum - \zA_2 \, \zB_1}_{T_3} 
+ \underbrace{\sum \zA_2 \, \zB_2}_{T_4}. \] 
Now 
\[ \begin{aligned} T_1 = &  
\sum\limits_{q=1}^{e-1}\binom{e-1}{q} (-1)^{e-1-q} \, b_1^{e-1-q} \, b_2^q \, b_\ux^{-1}
\, g_1^q \, g_2^{e-1-q} \, f_\ux^d \, g_\ux\\ 
= \, & (-1)^{e-1} \, b_1^{e-1} \, b_\ux^{-1} \, g_2^{e-1} \, f_\ux^d \, g_\ux 
\left\{ \lp 1-\frac{b_2 \, g_1}{b_1 \, g_2}\rp^{e-1}-1 \right\}\\
= \, & (-1)^{e-1} \, (b \, g)^{e-1} \, b_\ux^{-1} \, f_\ux^d \, g_\ux + 
(-1)^e \, b_1^{e-1} \, b_\ux^{-1} \, g_2^{e-1} \, f_\ux^d \, g_\ux, 
\end{aligned} \] 
and after similar calculations,
\[ \begin{aligned} T_2 = & 
(-1)^e \, b_1^{e-1} \, b_\ux^{-1} \, x_2^{-e+1} \, f_\ux^d \, g_\ux^e
+ (-1)^{e-1} \, b_1^{e-1} \, b_\ux^{-1} \, g_2^{e-1} \, f_\ux^d \, g_\ux \, , \\ 
T_3 = & 
(-1)^e \, (b \, f)^{e-1} \, b_\ux^{-1} \, f_\ux^{d-e+1} \, g_\ux^e+ 
(-1)^{e-1} \, b_1^{e-1} \, b_\ux^{-1} \, f_2^{e-1} \, f_\ux^{d-e+1} \, g_\ux^e \, , \\ 
T_4 = & (-1)^{e-1} \, b_1^{e-1} \, b_\ux^{-1} \, x_2^{-e+1} \, f_\ux^d \, g_\ux^e
+ (-1)^e \, b_1^{e-1} \, b_\ux^{-1} \, f_2^{e-1} \, f_\ux^{d-e+1} \, g_\ux^e. \\ 
\end{aligned} \] 
Notice that six of the eight terms cancel in pairs, for instance, 
the first term of $T_2$ cancels with 
the first term of $T_4$. We are left with 
\[ \begin{aligned} 
\tPsi & = (-1)^{e-1} \, (b \, g)^{e-1} \, b_\ux^{-1} \, f_\ux^d \, g_\ux
+ (-1)^e \, (b \, f)^{e-1} \, b_\ux^{-1} \, f_\ux^{d-e+1} \, g_\ux^e, \\ 
 & = \frac{(-1)^{e-1} \, f_\ux^{d-e+1} \, g_\ux}{b_\ux} \, 
[ \, (b \, g)^{e-1}  \, f_\ux^{e-1} - (b \, f)^{e-1}  \, g_\ux^{e-1}]. 
\end{aligned} \] 
Rewrite $b_\ux$ using the Pl{\"u}cker syzygy 
$b_\ux \, (f \, g) = (b \, g) \, f_\ux - (b \, f) \, g_\ux$, and factor the numerator.  
This gives 
\[ \tPsi = (-1)^{e-1} \, (f \, g) \, \sum\limits_{i=1}^{e-1} \, 
(b \, f)^{i-1} \, (b \, g)^{e-i-1} \, f_\ux^{d-i} \, g_\ux^i \, . \] 
Now make a change of variable $(b_1,b_2) = (z_2,-z_1)$. Then 
$(b \, f) = b_1 f_2 - b_2 f_1 = f_\uz, (b \, g) = g_\uz$, and 
$\DD_F \, \DD_G \, \tPsi = (-1)^{e-1} \, \Mor(\uz,\ux)$. By formula~(\ref{FGr-GFr}), 
\[ F \, \tG_r - G \, \tF_r  =
\frac{(-1)^e}{(e-2)!} \, 
(\frac{\partial}{\partial z_2})^r (-\frac{\partial}{\partial z_1})^{e-2-r} \, 
\DD_F \, \DD_G \, \tPsi, \] 
which is the same as  $\theta_r$ by formula~(\ref{expression.theta.r}). 
This completes the proof of identity~(\ref{morley.identity}), and hence 
that of Proposition~\ref{prop.morley}. \qed 

\bibliographystyle{plain}

\centerline{---}
\vspace{2cm} 

\parbox{7.5cm}{\small 
\, Carlos D'Andrea \\
Departament d'{\`A}lgebra i Geometria \\ 
Facultat de Mat{\`e}matiques \\ 
Universitat de Barcelona \\ 
Gran Via de les Corts Catalanes, 585 \\
E-08007 Barcelona, Spain. \\ 
{\tt carlos@dandrea.name} }
\parbox{7cm}{\small 
\, Jaydeep V.~Chipalkatti \\ 
433 Machray Hall \\ 
Department of Mathematics \\ 
University of Manitoba \\ 
Winnipeg R3T 2N2, Canada. \\ 
{\tt chipalka@cc.umanitoba.ca}}

\vspace{1.3cm} 

\parbox{6cm}{\small 
{Abdelmalek Abdesselam} \\ 
LAGA, Institut Galil\'ee \\ CNRS UMR 7539\\
Universit{\'e} Paris XIII\\
99 Avenue J.B. Cl{\'e}ment\\
F93430 Villetaneuse, France. \\
{\tt abdessel@math.univ-paris13.fr} }
\end{document}